\documentclass[12pt,reqno]{amsart}
\usepackage{amsmath, amssymb, amsthm, mathtools}
\usepackage{url}
\usepackage[breaklinks]{hyperref}

\renewcommand{\Im}{\operatorname{Im}}

\renewcommand{\Re}{\operatorname{Re}}
\renewcommand{\Im}{\operatorname{Im}}

\newcommand{\e}{\epsilon}

\renewcommand{\(}{\left\(}
\renewcommand{\)}{\right\)}
\renewcommand{\[}{\left\[}
\renewcommand{\]}{\right\]}

\numberwithin{equation}{section}
 \theoremstyle{plain}
\newtheorem{theorem}{Theorem}[section]
\newtheorem{lemma}[theorem]{Lemma}
\newtheorem{remark}[]{Remark}
\newtheorem{definition}[theorem]{Definition}

\newtheorem{corollary}[theorem]{Corollary}

   \makeatletter
\def\proof{\@ifnextchar[{\@oproof}{\@nproof}}
\def\@oproof[#1][#2]{\trivlist\item[\hskip\labelsep\textit{#2 Proof of\
#1.}~]\ignorespaces}
\def\@nproof{\trivlist\item[\hskip\labelsep\textit{Proof.}~]\ignorespaces}

\makeatother

\begin{document}
\title[A Lambert series associated with Siegel cusp forms]{An asymptotic expansion for a  Lambert series associated to Siegel cusp forms of degree $n$} 

\author{Babita}
\address{Babita\\Department of Mathematical Sciences\\ Indian Institute of Technology(Banaras Hindu University), Varanasi \\
221005, Uttar Pradesh, India.} 
\email{babita.rs.mat19@itbhu.ac.in}

\author{Abhash Kumar Jha}
\address{Abhash Kumar Jha\\Department of Mathematical Sciences\\Indian Institute of Technology(Banaras Hindu University), Varanasi \\
221005, Uttar Pradesh, India.} 
\email{abhash.mat@iitbhu.ac.in}


 \author{Bibekananda Maji}
\address{Bibekananda Maji\\ Department of Mathematics \\
Indian Institute of Technology Indore \\
Indore, Simrol, Madhya Pradesh 453552, India.} 
\email{bibek10iitb@gmail.com, bibekanandamaji@iiti.ac.in}

\author{Manidipa Pal}
\address{Manidipa Pal\\ Department of Mathematics \\
National Institute of Technology Calicut \\
Kozhikode,  Kerala,  673601,  India. }
\email{manidipa@nitc.ac.in,  math.manidipa@gmail.com}


\maketitle

\begin{abstract}
Utilizing inverse Mellin transform of the symmetric square $L$-function attached to Ramanujan tau function, Hafner and Stopple proved  a conjecture of Zagier,  
which states that the constant term of the automorphic function $y^{12}|\Delta(z)|^2$ i.e.,   the Lambert series $y^{12}\sum_{n=1}^\infty \tau(n)^2 e^{-4 \pi n y}$ 
can be expressed in terms of the non-trivial zeros of the Riemann zeta function.  This study examines certain Lambert series associated to Siegel cusp forms of degree $n$ twisted by a character $\chi$ and observes a similar phenomenon.  

\end{abstract}

\section{Introduction}

The famous Ramanujan tau function $\tau(n)$ was introduced and extensively studied by Ramanujan in his seminal paper \cite{ram36}.  It is the $n$th Fourier coefficient of  the Ramanujan delta function $\Delta(z)$, which is a cusp form of weight $12$ for the full modular group $\text{SL}_{2}(\mathbb{Z})$.  It is defined by
		$$\Delta(z):= q \prod_{n=1}^\infty (1-q^n)^{24} = \sum_{n=1}^\infty \tau(n)q^n, ~~q=e^{2\pi i z},  z \in \mathcal{H}.$$
Based on numerical values of $\tau(n)$,  Ramanujan observed several properties (congruences,  recursive relations and an estimate) for $\tau(n)$ which were later proved by mathematicians over the time.  More specially,  Ramanujan conjectured the estimate $\tau(n)=O\left(n^{11/2+\epsilon}\right)$ which was proved by Deligne in 1974.
 In 1981,  Zagier \cite[p.~417]{Zag}, \cite[p.~271]{Zag92} found another interesting connection of Ramanujan tau function with the non-trivial zeros of the Riemann zeta function $\zeta(s)$.   Mainly,  he conjectured that the constant term $c_0(y)$ of the automorphic function $\mathcal{F}(z):=y^{12}|\Delta(z)|^2$,  that is, the following Lambert series 
		\begin{equation}\label{Lambert series_Zagier}
			c_0(y):=y^{12} \sum_{n=1}^{\infty} \tau^2(n) \exp({-4 \pi ny}), 
		\end{equation} has an oscillatory behaviour when $y \to 0^{+}$.
	 Moreover, the Mellin transform of $c_0(y)$ is analytic except for poles at $s=1$ and $s=\rho/2$,  which indicates that $c_0(y)$ has the following asymptotic behaviour:
		\begin{align}\label{asymptotic expansion}
			c_0(y) \sim C + \sum_{\rho} y^{1 - \frac{\rho}{2}} A_{\rho},
		\end{align} 
	    as $y \to 0^{+}$,
		where the sum over $\rho$ runs through the non-trivial zeros of $\zeta(s)$, $ C$ is the average value of $\mathcal F$, precisely $C= \frac{3}{\pi} \langle \Delta , \Delta \rangle$,  and $A_\rho$ are certain complex numbers.  If $\rho$ is an $n$-fold zero,  then $A_\rho$ will be replaced by a polynomial in $\text{log} (y)$ of degree $n-1$.  
Furthermore,  Zagier mentioned that using the values of the Ramanujan tau function and the above asymptotic expansion,  one can numerically evaluate non-trivial zeros of $\zeta(s)$.  Under the assumption of the Riemann hypothesis,  the asymptotic expansion \eqref{asymptotic expansion} and the oscillatory behaviour of $c_0(y)$, for $y \rightarrow 0^{+}$, 
		\begin{align}\label{oscillatory behaviour}
			c_0(y) \sim C + y^{\frac{3}{4}} \sum_{m=1}^{\infty} a_m \cos \left(\frac{1}{2} \gamma_m \log(y) + \phi_m \right),
		\end{align}
		where $\gamma_m = \text{Im}(\rho)$, $a_m$ and $\phi_m$ are real constants, have been proved by Hafner and Stopple \cite{HS} by examining a heat kernel associated with Ramanujan tau function. 
		A similar phenomenon for Hecke eigenforms  for the full modular group as well as for the congruence subgroups have been observed in \cite{CKM} and \cite{CJKM}.  Subsequently,  an asymptotic expansion of an infinite series associated to Hecke-Maass eigenforms has been studied by Banerjee and Chakraborty \cite{BC-19}.  Recently, this problem in the case of Hilbert modular forms has been examined by Agnihotri \cite{Agni}. 
		 We refer  to \cite{JMS-21,  JMS-22,  MSS-22} for similar problems inspired from the above conjecture of Zagier.  
	Siegel modular forms are natural generalization of elliptic modular forms to higher dimension which were first introduced by Siegel to study quadratic forms. There are several connections between Siegel modular forms and other types of modular forms, viz, elliptic modular forms of integral and half-integral weight, Jacobi forms, etc. One such connection is the Saito-Kurokawa lift, which is a mapping from the space of elliptic modular forms of integer weight $2k -2$ to the space of Siegel modular forms of weight $k$ and degree $2.$ 

Recently,  Juyal and the first three authors \cite{BJJM24} studied an asymptotic behaviour of a Lambert series associated to Siegel cusp forms of degree $2$ involving their Fourier-Jacobi coefficients. 
 In this paper,  we extend the results in \cite{BJJM24} for degree $>2$.   Mainly,  
 we study an
 asymptotic expansion of a Lambert series associated to Siegel cusp forms of degree $n$ involving their Fourier-Jacobi coefficients with a twist by a Dirichlet character $\chi$.  Quite interestingly,  we observe a connection with the non-trivial zeros of the Dirichlet $L$-function $L(s, \chi^2)$. 

We now give outline of our paper. 
 In the next section,  we provide necessary backgrounds which are essential to state the main result of this article.  Section \ref{Main results} contains statement of the main result and its applications.  In Section \ref{Prili},  we state a few results without proof which will be used in the proof of the main result.  Finally,  in Section \ref{Proofs},  we give proof of results that are mentioned in Section \ref{Main results}.  

\section{Preliminaries}
Let $\mathbb C$ and $\mathcal H$ denote the complex plane and complex upper half-plane, respectively. For two matrices $A$ and $B$ of appropriate sizes, we define $A[B]:={\overline B}^tAB,$ where $\overline B^t$ is the conjugate transpose of the matrix $B.$ For $n>1$, we define $J:=\left(\begin{array}{cc} 0_n & {I_n}\\{-I_n} & 0_n \end{array}\right),$ where $0_n$ and $I_n$ denote the $n\times n$ zero and identity matrix,  respectively. The Siegel upper half-plane $ \mathcal{H}_n$ of degree $n$ is defined as 
$$
\mathcal{H}_n := \{ Z  \in M_{n \times n}(\mathbb{C})~|~Z^t=Z, ~ \Im (Z)~ \text{is positive definite}  \}. 
$$
The full Siegel modular group  $\Gamma_n:=Sp_{2n}(\mathbb Z)$  of degree $n$ is defined as 
\begin{equation*}
\Gamma_n:=\left\{M \in M_{2n \times 2n}(\mathbb{Z})~ | ~J[M]=J\right\}.
\end{equation*} 
If we write the matrix $M\in \Gamma_n$ as a block decomposition $M=\begin{pmatrix}
A  & B\\
C  & D\\
\end{pmatrix}$, where  $ A,B,C,D \in M_{n \times n}(\mathbb{Z}),$ then we have
\begin{equation*}
\Gamma_n=\left\{\begin{pmatrix}
A  & B\\
C  & D\\
\end{pmatrix} \bigg| A,B,C,D \in M_{n \times n}(\mathbb{Z}),A^tC = C^tA,  B^t D= D^t B,A^tD-B^tC=I_n\right\}.
\end{equation*}
The full Siegel modular group $\Gamma_n$ acts on the Siegel upper half-plane $\mathcal{H}_n$ as follows: 
$$
\begin{pmatrix}
A & B\\
C & D\\
\end{pmatrix} \cdot Z=(AZ+B)(CZ+D)^{-1},~\begin{pmatrix}
A & B\\
C & D\\
\end{pmatrix}\in \Gamma_n,~Z\in\mathcal H_n.
$$
\begin{definition}
Let $k$ and $n>1$ be positive integers. A complex-valued holomorphic function $F: \mathcal{H}_n \rightarrow \mathbb{C}$ is said to be a Siegel modular form of weight $k$ and degree $n$ if it satisfies 
$$
F (M \cdot Z) = {\det} (CZ+D)^{k} F( Z),
$$
for all $M = \begin{pmatrix}
A  & B\\
C  & D\\
\end{pmatrix} \in \Gamma_n$,  $Z \in \mathcal{H}_n,$ and has a Fourier series expansion of the form 
\begin{equation}\label{summation}
F(Z)= \sum_{\substack{T\ge 0}} A_F(T) e^{2\pi i (tr(TZ))}, 
\end{equation}
where $T$ runs through positive semidefinite half-integral (i.e., $T=(t_{ij}),~t_{ii}, 2t_{ij}\in \mathbb Z$) $n\times n$ matrices.
Further, we say  $F$ is a Siegel cusp form if the summation in \eqref{summation} runs over positive definite half-integral matrices $T$. 

\end{definition}
The space of Siegel modular forms and Siegel cusp forms of weight $k$ and degree $n$ are denoted by $M_k(\Gamma_n)$ and $S_k(\Gamma_n),$ respectively.




\subsection{Fourier-Jacobi expansion} 
Let $Z\in\mathcal H_n$ and $T$ be a positive semidefinite half-integral $n\times n$ matrix. Write the following block decomposition of $Z$ and $T:$
\begin{equation}\label{values of Z,T}
Z = \begin{pmatrix}
Z_1  & w\\
\bar w^t  & z\\
\end{pmatrix} , ~~~~~~~~~~~~~~~~~~~~~~~
T = \begin{pmatrix}
T_1  & \frac{1}{2}\lambda\\
\frac{1}{2} \bar \lambda^t  & m\\
\end{pmatrix},
\end{equation}
where $Z_1$ and $T_1$ denotes the upper left $(n-1)\times (n-1)$ corners of the matrix $Z$ and $T,$ respectively.  
Let $F\in M_k(\Gamma_n)$ with the Fourier series expansion as in \eqref{summation}. Rewrite the Fourier series expansion of $F$ as follows:
\begin{eqnarray*}\label{f-j}\nonumber
F(Z) &=& \sum_{T \ge 0} A_F(T) \exp({2\pi i(tr(TZ))})\\
  & =& \sum\limits_{\substack{T = \begin{pmatrix}
T_1  & \frac{1}{2}\lambda\\
\frac{1}{2}\bar \lambda^t  & m\\
\end{pmatrix} \ge 0}} \!\!\!\!\!\!\!\!\!\!A_F  \left( \begin{pmatrix}
T_1  & \frac{1}{2}\lambda\\
\frac{1}{2}\bar \lambda^t  & m\\
\end{pmatrix} \right)  \exp({2\pi i ( tr(T_1Z_{1}) +\bar \lambda^tw + mz)})  \\
&=&  \sum\limits_{m=0}^{\infty} \phi_{m}(Z_1,  w) \exp({2\pi i m z }),
\end{eqnarray*}
where for each $m$, the function $\phi_{m}(Z_1,  w)$ is called the $m$th Fourier-Jacobi coefficient of $F$ and is defined as
\begin{align*}
\phi_{m}(Z_1,  w) := \sum_{\substack{T_1 \ge 0,\lambda \in \mathbb{Z}^{n-1}\\ }} A_F  \begin{pmatrix}
T_1  & \frac{1}{2}\lambda\\
\frac{1}{2}\bar \lambda^{t}  & m\\
\end{pmatrix}  \exp({2\pi i ( tr(T_1 Z_1) +\bar \lambda^{t} w)}).  
\end{align*}
For each $m,$ the Fourier-Jacobi coefficient $\phi_{m}$ is a Jacobi form of weight $k$, index $m$ and degree $n-1$. It is easy to check that if $F\in S_k(\Gamma_n),$ then its Fourier-Jacobi expansion is given by $F(Z)=\sum\limits_{m=1}^{\infty} \phi_{m}(Z_1,  w) \exp({2\pi i m z })$ and each of the Fourier-Jacobi coefficients are Jacobi cusp forms. For more details on the theory of Siegel modular forms, we refer to \cite{Andrianov, 1-2-3,  klingen, maass}.  

\subsection{Dirichlet series associated to Siegel cusp forms}  Let  $F, G\in S_k(\Gamma_n)$ be Siegel cusp forms with Fourier-Jacobi expansions given by
\begin{equation}\label{fj-exp}
F(Z)=\sum\limits_{m=1}^{\infty} \phi_{m}(Z_1,  w) \exp({2\pi i m z }),~~~~~G(Z)=\sum\limits_{m=1}^{\infty} \psi_{m}(Z_1,  w) \exp({2\pi i m z }).
\end{equation}
The Rankin-Selberg  type Dirichlet series  associated to $F$ and $G$ is defined as
\begin{align}\label{Rankin-Selberg}
D_{F,  G}(s):=  \sum_{m=1}^\infty  \frac{\langle \phi_{m},  \psi_{m}  \rangle}{m^s},~~~~~~~~~~~~
\end{align}
where $\langle \phi_m , \psi_m \rangle$ is the Petersson scalar product of Jacobi cusp forms $\phi_m$ and $\psi_m.$ 
For precise definition of $\langle \phi_m , \psi_m \rangle,$  we refer to \cite{Ziegeler} for more details.

The series \eqref{Rankin-Selberg} was introduced and studied by Kohnen and Skoruppa \cite{k-s} for Siegel cusp forms of degree two. Later, Krieg \cite{krieg} studied the analytic properties of the series \eqref{Rankin-Selberg} associated to Siegel cusp forms of degree $n>2.$ 
Next, we consider the twist of the series \eqref{Rankin-Selberg} by a Dirchlet character. Let $N$ be a positive integer and $\chi$ a Dirichlet character  modulo $N$. 
The twist of the series $D_{F,  G}(s)$ by $\chi$ is defined as
\begin{align}\label{Rankin-Selberg-twisted}
\mathbb D_{F,  G,\chi}(s):=  \sum_{m=1}^\infty \chi(m) \frac{\langle \phi_{m},  \psi_{m}  \rangle}{m^s}
\end{align}

Kohnen \cite{Kohnen-twisted} studied the analytic properties of the series \eqref{Rankin-Selberg-twisted} associated to Siegel cusp forms of degree two. Later, the result of Kohnen \cite{Kohnen-twisted} was extended by Krieg, Kohnen and Sengupta \cite{KKS} to Siegel cusp forms of degree $n>2.$ 
 
\begin{remark}\cite[p.~246,  Lemma 1]{krieg}
The Petersson scalar product of the Fourier-Jacobi coefficients $\phi_m$ and $\psi_m$ satisfies
$$ \langle \phi_m,\psi_m \rangle = O(m^k).$$
\end{remark}

In view of the above remark, we see that the series \eqref{Rankin-Selberg} and \eqref{Rankin-Selberg-twisted} are absolutely convergent for $\Re(s)>k+1$. 

{Further,  Kohnen \cite[p.~718]{Koh93},   \cite[p.~134]{Koh11} conjectured and Kohnen and Sengupta \cite{KS17} proved that for any $\epsilon>0$,  
\begin{align}\label{Ram-Pet}
\langle \phi_m,  \phi_m \rangle = O_{F,  \epsilon} \left( m^{k-1+\epsilon} \right).  
\end{align}
when $F\in S_k(\Gamma_2)$ is a Hecke eigenform which is a Saito-Kurokawa lift.  Recently,  Kumar and Paul \cite{KP21} proved that the above conjecture is true on average for arbitrary degree.  Assuming \eqref{Ram-Pet},  one can clearly see that the Dirichlet series $D_{F ,  G}(s)$ and $\mathbb D_{F , G}(s,\chi)$ are absolutely convergent for $\Re(s)>k$.}


\section{Statement of Main Results}\label{Main results}
In this section, we state our main results. 
The main aim of this paper is to study certain Lambert series associated to two Siegel cusp forms of degree $n$. More precisely,  let $\chi$ be a Dirichlet character modulo $N$ and $F,G \in S_k(\Gamma_n)$ with the Fourier-Jacobi  expansion as in \eqref{fj-exp}.  For $\alpha>0,$ we are interested to study the following Lambert series:
\begin{align}\label{Lambert_Siegel}
&\sum_{m=1}^\infty \chi(m) \langle \phi_m,\psi_m \rangle \exp(-4 \pi m \alpha),
\end{align}
where $\phi_m$ and $\psi_m$ are $m$th Fourier-Jacobi coefficients of $F$ and $G$, respectively.
Before stating the main result of this paper,  we first define an arithmetical function $a_{F,  G}(m,\chi)$ that satisfy the following generating function,  for $\Re(s)>k$,  
\begin{align}\label{series1}
 \sum_{m=1}^\infty \frac{a_{F,  G}(m,\chi)}{m^s}= \frac{L(2s-2k+2n,\bar \chi^2) \mathbb D_{F, G, \bar \chi}(s)}{L(2s-2k+1,\bar \chi^2)},
\end{align}
 where $L(s,\chi)=\sum_{m=1}^{\infty}\chi(m)m^{-s}$ stands for Dirichlet $L$-function and $\mathbb{D}_{F, G,  \chi}(s)$ is the twisted Rankin-Selberg $L$-function defined in \eqref{Rankin-Selberg-twisted}.  
Note that for $N=1$, $\bar \chi(t) = \bar \chi^2(t)=1, \mathbb D_{F,G,\bar\chi}(s)=D_{F,G}(s)$ and $L(s,\bar \chi^2)=\zeta(s).$

 For complex numbers $\kappa$ and $\mu$,   let us consider the second order ordinary differential equation:
\begin{align*}
\frac{d^2w}{dz^2} + \left( \frac{1-4\mu^2+4 \kappa z -z^2}{4 z^2}  \right)w=0. 
\end{align*}
One of the solutions of the above differential equation is the Whittaker function $W_{\kappa,  \mu}(z)$.  

With the above notation and definitions in our hand,  we are ready to state the main result of this article.

\begin{theorem}\label{main theorem_1}
Let $F$ and $G$ be Siegel cusp forms of weight $k$ and degree $n$  with the Fourier-Jacobi expansion given  as in \eqref{fj-exp} and satisfying \eqref{Ram-Pet}. Let $\alpha$ and $\beta$ be two positive real numbers such that $\alpha \beta =1$. 
~Let $\chi$ be a Dirichlet character modulo $N $ such that $\chi^2$ is primitive. Under the assumption of the simplicity hypothesis of the non-trivial zeros of $L(s,  \chi^2)$,  we have
\begin{align*}
&\sum_{m=1}^\infty \chi(m)\langle \phi_m,\psi_m \rangle \exp(-4 \pi m \alpha) \\
& =   \left(\frac{\beta^{2k-n}N^{2n-2k+2}}{ \pi^{n-\frac{1}{2}} g({\bar \chi})^4 g({\chi^2})}\right) \sum_{m=1}^\infty a_{F,  G}(m,\chi) \left(\frac{4 \pi m \beta}{N^2} \right)^{\frac{n-1-k}{2}} W_{\frac{n+k}{2}-1,  \frac{k-n}{2}}\left(\frac{4\pi m \beta}{N^2} \right) \exp \left(\frac{-2\pi m \beta}{N^2} \right) \nonumber \\
&  + \sum_{\rho}\ \frac{ \mathbb D^*_{F, G,\chi}\left(\frac{\rho}{2}+k-n \right) (4 \pi \alpha)^{n-k-\frac{\rho}{2} } }{(2\pi / N)^{2n-2k-\rho} \Gamma(\frac{\rho}{2}) L'(\rho,\chi^2)} + \mathcal{R}_k,  \nonumber 
\end{align*}
where the term $\mathcal {R}_{k}$ is defined as
\begin{align}
\mathcal {R}_{k}= \begin{cases} \frac{(-1)^{n+1}(2n)! \langle F , G \rangle}{(4\pi)^n \alpha^k (n-1)! B_{2n}},  & ~\text{if}~ N=1,  \\
0,  & otherwise,  
\end{cases}
\end{align}
and $L'(s,\chi^2)$ denotes the derivative of $L(s,\chi^2)$ and the sum over $\rho$ runs through non-trivial zeros of $L(s,\chi^2)$ satisfying the bracketing condition,  that is, 
the terms corresponding to $\rho_1$ and $\rho_2$ are included in the same bracket if they satisfy
\begin{align} \label{bracketing condition for L}
|\Im(\rho_1) - \Im(\rho_2)| < \exp \left( -\frac{C_0 |\Im(\rho_1)|}{\log(|\Im(\rho_1)|+3 )} \right) + \exp \left( -\frac{C_0 |\Im(\rho_2)|}{\log(|\Im(\rho_2)|+3 )} \right),
\end{align}
where $C_0$ is some positive constant.
\end{theorem}
In particular,  if we consider $\chi$ to be the trivial character with modulus $N=1$.  Then we have the following result.

\begin{corollary}\label{main theorem_2}
Let $F,  G,  \alpha,  \beta$ and $\chi$ be as in   Theorem \ref{main theorem_1}. 
 Under the assumption of the simplicity hypothesis of the non-trivial zeros of $\zeta(s)$,  we have
\begin{align}\label{main eqn}
&\sum_{m=1}^\infty \langle \phi_m,\psi_m \rangle \exp(-4 \pi m \alpha) \nonumber \\
& =   \frac{ \beta^{2k-n}}{\pi^{n-\frac{1}{2}}} \sum_{m=1}^\infty a_{F,  G}(m) (4 \pi m \beta)^{\frac{n-1-k}{2}} W_{\frac{n+k}{2}-1,  \frac{k-n}{2}}(4\pi m \beta) \exp(-2\pi m \beta) \nonumber \\
& + \frac{(-1)^{n+1}(2n)! \langle F,G \rangle}{(4\pi)^n \alpha^k (n-1)! B_{2n}} + \sum_{\rho}\frac{ D^*_{F, G}\left(\frac{\rho}{2}+k-n \right) (4 \pi \alpha)^{n-k-\frac{\rho}{2} } }{(2\pi)^{2n-2k-\rho} \Gamma(\frac{\rho}{2}) \zeta'(\rho)},  
\end{align}
where the sum over $\rho$ runs through non-trivial zeros of $\zeta(s)$ satisfying the bracketing condition \eqref{bracketing condition for L} and 
$B_{2n}$ denotes the $2n$th Bernoulli number.
\end{corollary}

As an application of the above corollary,  we have the following asymptotic expansion of the Lambert series.  

\begin{corollary}\label{application_asymptotic}
Let $F$ and $G$ be as in Theorem \ref{main theorem_1}. 
If $ \langle F,  G \rangle \neq 0$,  then for $\alpha \rightarrow 0^{+} $,  we have
\begin{align}\label{asymptotic}
\alpha^k \sum_{m=1}^\infty \langle \phi_m,\psi_m \rangle \exp(-4 \pi m \alpha) \sim \frac{(-1)^{n+1}(2n)! \langle F,G \rangle}{(4\pi)^n  (n-1)! B_{2n}}. 
\end{align}
\end{corollary}


\section{Some Well-known Results}\label{Prili}
In this section,  we state a few results without proof which will play a crucial role in the proof of the Theorem {\rm \ref{main theorem_1}. 

Let $\chi$ be a Dirichlet character modulo $N$  and $F,G\in S_k(\Gamma_n).$  
Consider the completed Dirichlet series $D^*_{F,G}(s)$ and $\mathbb D^*_{F,  G,\chi}(s)$ associated to $F$ and $G$ defined by
\begin{align}
&D^{*}_{F,  G}(s):= (2\pi)^{-2s}\Gamma(s)\Gamma(s+n-k)\zeta(2s+2n-2k)D_{F,G}(s),  \label{complete_eqn} \\
&\mathbb D^{*}_{F,  G,\chi}(s) := \left(\frac{2\pi}{N}\right)^{-2 s} \Gamma(s) \Gamma(s+n-k) L(2s+2n-2k,\chi^2) \mathbb D_{F,  G,\chi}(s).  \label{comlete_eqn_character}
\end{align}


Krieg \cite{krieg} studied the analytic properties of $D^{*}_{F,  G}(s)$ and proved the following result. 

\begin{theorem} \cite[p.~249]{krieg}\label{Functional equn}
Let $F$ and $G $ be Siegel cusp forms of weight $k$ and degree $n.$ Then $D^{*}_{F,  G}(s)$ is a holomorphic function of $s \in \mathbb C$ except for possible simple poles at $s=k$ and $s=k-n$ with residues $\frac{1}{2}\pi^{n-k}\langle F,G \rangle$ and $-\frac{1}{2}\pi^{n-k}\langle F,G \rangle,$ respectively and satisfies the functional equation $D^{*}_{F,  G}(2k-n-s)=D^{*}_{F,  G}(s)$.
\end{theorem}

Later,   Kohnen,  Krieg and Sengupta \cite{KKS} extended the above result by studying the analytic properties of $\mathbb D^{*}_{F,  G,\chi}(s) $.    More precisely,   they proved the following result.  

\begin{theorem} \cite[p.~495]{KKS}\label{residue}
Let $F$ and $G $ be Siegel cusp forms of weight $k$ and degree $n,$ and $\chi$ be a Dirichlet character modulo $N>1.$ If $\chi$ is not principal, then $\mathbb D^*_{F,  G,\chi}(s)$ extends to an entire function of $s$.
Further, if  $\chi^2$ is primitive, then $\mathbb D^*_{F,  G,\chi}(s)$ satisfies the following functional equation:
 $$\mathbb D^{*}_{F,  G,\chi}(2k-n-s)=\left(\frac{g(\chi)}{\sqrt N} \right)^4 \mathbb D^{*}_{F,  G, \bar \chi}(s),$$
  where $g(\chi):= \sum\limits_{v \pmod{N}} \!\!\!\!\chi(v) e^{2 \pi i v / N}$ denotes the Gauss sum associated to the character $\chi.$ 
\end{theorem}

  


The next theorem gives analytic properties of Dirichlet $L$-functions.
\begin{theorem} \cite[p.~84]{IK}\label{L-Functional equn_2}
Let $\chi$ be a primitive Dirichlet character modulo $N \ge 1.$ Then $L(s,\chi)$  extends to a meromorphic function on $\mathbb C$ which is entire if $\chi \ne 1$ and otherwise admits a unique simple pole with residue $1$ at $s=1$. The completed $L$-function
  \begin{align*}
\Lambda(s,\chi)=\left (\frac{N}{\pi} \right)^{s/2} \Gamma \left( \frac{s+\e}{2} \right) L(s,\chi)
  \end{align*}
  is entire if $N \ne 1$ and has simple poles with residue $1$ at $s=0$ and $s=1$ otherwise. Moreover, it satisfies the following functional equation
   \begin{align*}
   \Lambda(s,\chi) = i^{-\e}~ \frac{g(\chi) }{\sqrt N} \Lambda(1-s,\bar \chi), ~~~~~~~~~~~~where~\e =\frac{1}{2} (1-\chi(-1)).
    \end{align*}  

\end{theorem}

 Next,  we define a special function known as Meijer $G$-function.
Let $m,n,p,q$ be integers with $0\leq m \leq q$, $0\leq n \leq p$.  Let $a_1, \cdots, a_p$ and $b_1, \cdots, b_q$ be $p+q$ complex numbers such that $a_i - b_j \not\in \mathbb{N}$, for $i \in [1, n]$,  $ j \in [1,  m]$.  The Meijer $G$-function \cite[p.~415]{NIST} is defined by the following line integral:
\begin{align}\label{MeijerG}
	G_{p,q}^{\,m,n} \!\left(  \,\begin{matrix} a_1,\cdots , a_p \\ b_1, \cdots , b_q \end{matrix} \; \Big| z \right) = \frac{1}{2 \pi i} \int_L \frac{\prod_{j=1}^m \Gamma(b_j + w) \prod_{j=1}^n \Gamma(1 - a_j -w) z^{-w}  } {\prod_{j=m+1}^q \Gamma(1 - b_j - w) \prod_{j=n+1}^p \Gamma(a_j + w)}\mathrm{d}w.  
\end{align}
Here, we assume that the line of integration $L$ separates the poles of the factors $\Gamma(b_j+w)$ from the poles of the factors $\Gamma(1-a_j-w)$. The above integral converges absolutely if $2(m+n)> p+q$
 and $ |\arg(z)| < (2m+2n - p-q) \frac{\pi}{2}$.
 Next, we have the following special evaluation of the Meijer $G$-function.  
 
\begin{lemma}\cite[p.~58]{KT} \label{Special case of Meijer G}
 For $|\arg(z)|< \frac{\pi}{2}$,  we have
\begin{align}
G_{1,2}^{\,2,0} \!\left(  \,\begin{matrix} a_1 \\ b_1,  b_2 \end{matrix} \; \Big| z \right) = z^{\frac{b_1+b_2-1}{2}} e^{-z/2} W_{\kappa, \mu}(z),
\end{align}
where $W_{\kappa, \mu}(z)$ is the Whittaker function with $\kappa= \frac{1}{2}(b_1+b_2-1)-a_1$ and $\mu=\frac{1}{2}(b_1 - b_2)$.
\end{lemma}

The next lemmas tells us about the behaviour of the Whittaker function.  
\begin{lemma} \cite[p.~69,  eq. (13.14.21)]{NIST}  \label{Asym_Whittaker}
As $z \rightarrow \infty$ with $|\arg(z)| <  \frac{3\pi}{2}$,  one has

\begin{align*}
W_{\kappa, \mu}(z) \sim \exp(-z/2) z^{\kappa}.
\end{align*}
\end{lemma}



The next lemma gives Stirling's bound for $\Gamma(z)$.
\begin{lemma} \cite[p.~151,~A.4.]{IK}   \label{Stirling}
	Let $z = \sigma + i {T}$ with $p \leq \sigma \leq q.$ Then  we have
	\begin{equation}\label{Stirling_equn}
		|\Gamma (\sigma + i {T})| \sim \sqrt{2\pi} | {T}|^{\sigma - 1/2} e^{-\frac{1}{2} \pi |{T}|},   \quad {\rm as} \quad |{T}|\rightarrow \infty.
	\end{equation}
\end{lemma}

 \begin{lemma}\cite[p. 8]{MSS-22}\label{bound_for_L}
Assume that there exist a sequence of positive numbers $T_n$ with arbitrary large absolute value satisfying $|T_n-\Im(\rho)|>\exp(-A|\Im(\rho)|)/\log(|\Im(\rho)|)$ for every non trivial zero $\rho$ of $L(s,\chi)$, where $A$ is some suitable positive constant. Then for $0<B<\frac{\pi}{4},$ we have
 \begin{align*}
 \frac{1}{L(\sigma+iT_n,\chi)} < e^{BT}.
 \end{align*}
\end{lemma}
In the next section,  we present proof of all the main results of this paper.

\section{Proof of main results}\label{Proofs}

\begin{proof}[Theorem {\rm \ref{main theorem_1}}][]
Using the definition of $\Gamma(s)$ and $\mathbb D_{F,  G,\chi}(s)$,  it can be shown that,  $\Re(s)>k+1$,  
\begin{align}
\Gamma(s) \mathbb D_{F,  G,\chi}(s) =  \int_{0}^\infty \sum\limits_{m=1}^\infty \chi(m) \langle \phi_m,\psi_m \rangle \exp(- m x) x^{s-1} \mathrm{d}x.  
\end{align}
Thus, by inverse Mellin transform,  for any real number $c>k+1$,  one has
\begin{align}\label{right_vertical_integral_2}
&\sum_{m=1}^\infty \chi(m) \langle \phi_m,\psi_m \rangle \exp(-4 \pi m \alpha)
  =\frac{1}{2\pi i}\int_{c-i \infty}^{c+ i \infty}\Gamma(s)  \mathbb D_{F,  G,\chi}(s) (4 \pi \alpha)^{-s}  \mathrm{d}s \nonumber \\
 & = \frac{1}{2\pi i}\int_{(c)}\frac{\mathbb D^*_{F,  G,\chi}(s) (4 \pi \alpha)^{-s}}{(2\pi /N)^{-2s}\Gamma(s+n-k) L(2s+2n-2k,\chi^2)}  \mathrm{d}s,
\end{align}
where the symbol $(c)$ denotes the line integral $c- i \infty$ to $c+ i \infty$.  
In view of Theorem \ref{residue}, 
 it is easy to observe that in \eqref{right_vertical_integral_2},  the trivial zeros of  $L(2s+2n-2k,\chi^2)$ are neutralised by the poles of $\Gamma(s+n-k)$.  Nevertheless, the non-trivial of zeros of $L(2s-2k+2n,\chi^2)$ will give us infinitely many poles of the integrand in the strip $k-n < \Re(s) < k -n+ \frac{1}{2}$.  Thus,  we construct a closed rectangular contour $\mathcal{C}$ consisting of the end points $c-i T,  c+i T,  c_1 + i T,  c_1 - iT$,  where $T$ is some large positive constant.  Here, we consider $c_1 \in (k-n-1,  k-n)$ so that the non-trivial zeros of $L(2s-2k+2n,\chi^2)$ lie inside  the contour $\mathcal{C}$.
 
 When $N=1$,  the Dirichlet character $\chi$ becomes trivial and Dirichlet $L$-function reduces to the Riemann zeta function.  In view of Theorem \ref{Functional equn},  one can check that the integrand function will have an extra pole at $s=k$.  However,  there would not be a pole at $s=k-n$. 
 
 Now,  apply the residue theorem to obtain 
\begin{equation}\label{CRT_2}
\frac{1}{2\pi i} \int_{\mathcal{C}} \frac{(2\pi / N)^{2s} \mathbb D^*_{F,G,\chi}(s)(4 \pi \alpha)^{-s}}{\Gamma(s+n-k) L(2s-2k+2n,\chi^2)}   \mathrm{d}s  = \sum_{ |\Im(\rho)| < T} \mathcal R_{\rho} + \mathcal R_{k}, 
\end{equation}
where the term $\mathcal R_{\rho}$ denotes the residual term at the non-trivial zero $\rho$ of $L(s,\chi^2)$. 
The residual term $\mathcal {R}_k$ is given by
\begin{align}\label{residue_k}
\mathcal {R}_{k}= \begin{cases} \frac{(-1)^{n+1}(2n)! \langle F , G \rangle}{(4\pi)^n \alpha^k (n-1)! B_{2n}},  & ~\text{if}~ N=1,  \\
0,  & otherwise.  
\end{cases}
\end{align}
To evaluate the above residual term we have used the fact that $\zeta(2n)=\frac{(-1)^{n+1}B_{2n}(2\pi)^{2n}}{2(2n)!}$, where $B_{2n}$ is the $2n$-th Bernoulli number.  
Under the assumption that the non-trivial zeros of Dirichlet $L$-functions are simple, we evaluate $\mathcal R_\rho$ as follows:
\begin{align}\label{residue_non-trivial_2}
\mathcal R_{\rho}= \frac{\mathbb D^*_{F, G,\chi}\left(\frac{\rho}{2}+k-n \right) (4 \pi \alpha)^{n-k-\frac{\rho}{2} } }{(2\pi / N)^{2n-2k-\rho} \Gamma\left(\frac{\rho}{2}\right) L'(\rho,\chi^2)},
\end{align}
where $L'$ denotes the derivative of the Dirichlet $L$-function.
Employing Stirling's bound \eqref{Stirling_equn} and Lemma \ref{bound_for_L},  horizontal integrals will tend to zero as $T \rightarrow \infty$.     Therefore,  making use of \eqref{right_vertical_integral_2} in \eqref{CRT_2} it yields that
\begin{align}\label{after CRT_2}
&\sum_{m=1}^\infty \chi(m) \langle \phi_m,\psi_m \rangle \exp(-4 \pi m \alpha) \nonumber \\
 & =\frac{1}{2\pi i}\int_{(c_1)}\frac{\mathbb D^*_{F,  G,\chi}(s) (4 \pi \alpha)^{-s}}{(2\pi /N)^{-2s}\Gamma(s+n-k) L(2s-2k+2n,\chi^2)}  \mathrm{d}s + \sum_{\rho} \mathcal R_{\rho} + \mathcal {R}_{k}.
\end{align}
Note that the sum over $\rho$ runs through all the non-trivial zeros of the Dirichlet L-function $L(s, \chi^2)$.  The convergence of this kind of infinite sum is not known without bracketing condition \eqref{bracketing condition for L}.  A similar infinite series appeared in the work of Hardy and Littlewood \cite[p.~156]{HL-1916}.   To know more about this kind of sums,  readers can see \cite{AGM22,  GM23}. 

Next, we concentrate on the left vertical integral, 
\begin{align}\label{left_vertical_integral_2}
\mathbb V_k(\alpha,\chi):&= \frac{1}{2\pi i}\int_{(c_1)}\frac{\mathbb D^*_{F,  G,\chi}(s) (4 \pi \alpha)^{-s}}{(2\pi /N)^{-2s}\Gamma(s+n-k)L(2s-2k+2n,\chi^2)}  \mathrm{d}s \nonumber \\
&=\frac{1}{2\pi i}\int_{(c_1)}\frac{(2\pi /N)^{2s} \mathbb D^*_{F,  G,\bar \chi}(2k-n-s) (4 \pi \alpha)^{-s}}{\left(\frac{g(\bar \chi)}{\sqrt N}\right)^4 \Gamma(s+n-k)L(2s-2k+2n,\chi^2)}  \mathrm{d}s.
\end{align}
here we utilized Theorem \ref{residue}.  Further,  use \eqref{comlete_eqn_character} to see that
\begin{align}\label{Rankin_functional_2}
& \mathbb D^*_{F,  G,\bar \chi}(2k-n-s) \nonumber \\
&=(2 \pi /N)^{-2(2k-n-s) } \Gamma(2 k-s-n) \Gamma(k-s) L(2k-2s,\bar \chi^2) \mathbb D_{F,  G,\bar \chi}(2 k-n-s).
\end{align}
Employ \eqref{Rankin_functional_2} in \eqref{left_vertical_integral_2} to deduce that
\begin{align}\label{v-k-alpha}
\mathbb V_k(\alpha,\chi)= \frac{1}{2\pi i}\int_{(c_1)}  \frac{\Gamma(2 k-n-s) \Gamma(k-s) L(2k-2s,\bar \chi^2) \mathbb D_{F,  G,\bar \chi}(2 k-n-s)(4 \pi \alpha)^{-s}}{ \left(\frac{2\pi}{N}\right)^{4k-4s-2n}\left(\frac{ g(\bar \chi)}{\sqrt N}\right)^4 \Gamma\left( s-k+n  \right) L(2s -2k+2n,\chi^2)} \mathrm{d}s.
\end{align} 
We change the variable by letting $2k-n-s=w$ in \eqref{v-k-alpha} to simply further.  Thus, 
\begin{align}\label{change variable_2}
&\mathbb V_k(\alpha,\chi) \nonumber\\
&= \left(\frac{\pi}{N^2 \alpha}\right)^{2k-n}   \frac{1}{2\pi i}\int_{(d_1)} \frac{\Gamma(w) \Gamma(w-k+n) L(2w+2n-2k,\bar \chi^2) \mathbb D_{F, G,\bar \chi}(w)}{\left(\frac{g(\bar \chi)}{\sqrt N}\right)^4 \Gamma\left(k-w  \right) L(2k-2w,\chi^2)}\left(\frac{\alpha N^4}{4\pi^3}\right)^w \mathrm{d}w, 
\end{align}
where  $k < \Re(w)=d_1 <k+1$. 
Use functional equation of the Dirichlet $L$-function,  i.e.,  Theorem \ref{L-Functional equn_2} to see that
\begin{align}\label{functional_eqn_after_substitution}
L(2k-2w,\chi^2)=\frac{g(\chi^2)}{\sqrt N} \frac{(N/ \pi)^{\frac{1}{2}-2k+2w}\Gamma \left(w-k+\frac{1}{2} \right)L(2w-2k+1,\bar \chi^2)}{\Gamma(k-w)}.
\end{align}
Substituting \eqref{functional_eqn_after_substitution} in \eqref{change variable_2},  we derive that
\begin{align}
&\mathbb V_k(\alpha,\chi)= \left(\frac{\alpha^{n-2k}N^{2n-2k+2}}{\pi^{n-\frac{1}{2}} g(\bar \chi)^4 g(\chi^2)}\right) \nonumber  \\
& \times  \frac{1}{2\pi i}\int_{(d_1)} \frac{\Gamma(w) \Gamma(w-k+n) L(2w+2n-2k,\bar \chi^2) \mathbb D_{F, G,\bar \chi}(w)}{\Gamma \left(w-k+\frac{1}{2}\right) L(2w-2k+1,\bar \chi^2)}\left(\frac{4 \pi }{\alpha N^2}\right)^{-w} \mathrm{d}w.  \label{invoking functional equation_2}
\end{align}
For $\Re(w)>k$,  we know that the Dirichlet series $\mathbb D_{F, G, \bar \chi}(w)$,  $L(2w+2n-2k,\bar \chi^2)$ and $1/ L(2w-2k+1,\bar \chi^2)$ are absolutely and uniformly convergent.  Thus,  one can write 
\begin{align}\label{series_2}
L(2w+2n-2k,\bar \chi^2) \frac{\mathbb D_{F, G,\bar \chi}(w)}{L(2w-2k+1,\bar \chi^2)} = \sum_{m=1}^\infty \frac{a_{F,  G}(m, \chi)}{m^w},~~\Re(w)>k,
\end{align}
where $a_{F,G}(m, \chi)$ is defined in \eqref{series1}.  Now,  substituting \eqref{series_2} in \eqref{invoking functional equation_2} and upon simplification,  for $\alpha \beta =1$,  we arrive
\begin{align}\label{Final_V_k_2}
\mathbb V_k(\alpha,\chi) = \left(\frac{\beta^{2k-n}N^{2n-2k+2}}{ \pi^{n-\frac{1}{2}} g(\bar \chi)^4 g(\chi^2)}\right) \sum_{m=1}^\infty a_{F,  G}(m,\chi) \mathbb I_{k}(m,  \beta),
\end{align}
where
\begin{align*}
\mathbb  I_{k}(m,  \beta):= \frac{1}{2\pi i}\int_{(d_1)} \frac{\Gamma(w) \Gamma(w-k+n) }{\Gamma\left(w-k+\frac{1}{2}  \right)}\left(\frac{4m\pi \beta}{N^2} \right)^{-w} \mathrm{d}w,
\end{align*}
Recalling the definition \eqref{MeijerG},  one can show that
 the above integral is nothing but
\begin{align}\label{I_k}
\mathbb I_{k}(m,  \beta)= G_{1,2}^{2,0}\left(\begin{array}{c}
\frac{1}{2}-k \\
0,n-k
\end{array} \Big| \frac{4 \pi m \beta}{N^2} \right). 
\end{align}
Utilize the relation between Meijer $G$-function and Whitteker function,  i.e.,  Lemma \ref{Special case of Meijer G} to see that
\begin{align}\label{Final_I_k_2}
\mathbb I_{k}(m,  \beta)= \left (\frac{4 \pi m \beta}{N^2}\right)^{\frac{n-1-k}{2}} W_{\frac{n+k}{2}-1,  \frac{k-n}{2}}\left(\frac{4\pi m \beta}{N^2} \right) \exp \left(\frac{-2\pi m \beta}{N^2} \right).
\end{align}
Finally,  putting \eqref{Final_I_k_2} in \eqref{Final_V_k_2}, the left vertical integral simplifies to 
\begin{align}\label{Final expression_V_k_2}
&\mathbb V_{k}(\alpha,\chi) \nonumber \\
& =  \frac{ \sqrt \pi \beta^{2k-n}N^{2n-2k+2}}{\pi^{n}g(\bar \chi)^4 g(\chi^2)} \!\!\sum_{m=1}^\infty \!\!\!a_{F,  G}(m,\chi) \!\left (\frac{4 \pi m \beta}{N^2}\right)^{\frac{n-1-k}{2}} \!\!\!\!\!\!W_{\frac{n+k}{2}-1,  \frac{k-n}{2}}\left(\frac{4\pi m \beta}{N^2} \right) \exp \left(\frac{-2\pi m \beta}{N^2} \right).
\end{align}
The convergence of the above series is an easy consequence of the fact that $a_{F,  G}(m,\chi)$ has a polynomial growth and the Whittaker function decays exponentially (see Lemma \ref{Asym_Whittaker}).  Finally,  the proof of  of Theorem \ref{main theorem_1} follows by making use  \eqref{Final expression_V_k_2} of the left vertical integral $V_{k}(\alpha,\chi)$  in \eqref{after CRT_2}, together with the residual terms \eqref{residue_k} and \eqref{residue_non-trivial_2}.   

\end{proof}

\begin{proof}[Corollary {\rm \ref{main theorem_2}}][]
The proof of this result immediately follows by letting $N=1$ in Theorem \ref{main theorem_1}. 
\end{proof}

\begin{proof}[Corollary {\rm \ref{application_asymptotic}}][]
Multiply  $\alpha^k$ on the both sides of \eqref{main eqn} and use $\alpha \beta=1$ to see that
\begin{align}\label{coro eqn}
&\alpha^k\sum_{m=1}^\infty \langle \phi_m,\psi_m \rangle \exp(-4 \pi m \alpha) \nonumber \\
& =   \frac{ \beta^{k-n}}{\pi^{n-\frac{1}{2}}} \sum_{m=1}^\infty a_{F,  G}(m) (4 \pi m \beta)^{\frac{n-1-k}{2}} W_{\frac{n+k}{2}-1,  \frac{k-n}{2}}(4\pi m \beta) \exp(-2\pi m \beta) \nonumber \\
& + \frac{(-1)^{n+1}(2n)! \langle F,G \rangle}{(4\pi)^n (n-1)! B_{2n}} + \sum_{\rho}\frac{ D^*_{F, G}\left(\frac{\rho}{2}+k-n \right) (4 \pi)^{n-k-\frac{\rho}{2}} \alpha^{n-\frac{\rho}{2} }} {(2\pi)^{2n-2k-\rho} \Gamma(\frac{\rho}{2}) \zeta'(\rho)}.  
\end{align}
Now,  considering the first infinite series on the right side of \eqref{coro eqn} and using Lemma \ref{Asym_Whittaker},  for $\alpha \rightarrow 0^{+}$, that is,  $\beta \rightarrow \infty$,  one can check that
\begin{align}\label{coro eqn 1}
 & \frac{ \beta^{k-n}}{\pi^{n-\frac{1}{2}}} \sum_{m=1}^\infty \frac{ a_{F,  G}(m)}{ (4 \pi m \beta)^{\frac{k-n+1}{2}}}  \exp(-2\pi m \beta)  W_{\frac{n+k}{2}-1, \frac{k-n}{2}}(4\pi m \beta) \nonumber \\
  & \ll \beta^{k-n} \sum_{m=1}^\infty \frac{ a_{F,  G}(m)}{ (4 \pi m\beta)^{\frac{k-n+1}{2}}}  \exp(-4\pi m \beta)  (4\pi m \beta)^{\frac{n+k}{2}-1}  \ll \frac{1}{\beta^K}, 
\end{align}
where $K>0$ is a large number.    Thus,  we can clearly see that the above sum tends to zero.  Next,  we can see that the infinite sum over $\rho$  on the right side of \eqref{coro eqn} will also vanish.  As
we know that the non-trivial zeros $\rho$ of $\zeta(s)$ satisfy $0<\Re(\rho)<1$.   Thus,  one has $n-\frac{1}{2} < \Re\left(n- \frac{\rho}{2}\right) <n$ and hence the  infinite sum over $\rho$ vanishes as $\alpha \rightarrow 0^{+}$.  This completes the proof of this corollary. 
\end{proof}

\begin{section}{Acknowledgement}
The second author is partially supported by SERB Start-up Research Grant (File No. SRG/2022/000487)  and MATRICS grant (File No. MTR/2022/000659).  The last author's research  supported by SERB MATRICS grant (File No. MTR/2022/000545) and SERB CRG grant (File No.  CRG/CRG/2023/002122). 
  Both the authors sincerely thank SERB for the support. 
\end{section}

\end{document}